\documentclass{article}
\usepackage{graphicx} 
\usepackage{amsmath, amstext, amsthm, amssymb, latexsym}
\usepackage[a4paper, margin=1in]{geometry} 
\usepackage{multicol}

\title{\textbf{Generalization of Arithmetico-Geometric Series and the Expectation Value of a $k$-Run of a Bernoulli Trial}}
\author{Priyansh Khare}
\date{\small 41, Sriram Colony, Hoshangabad Road, Bhopal, Madhya Pradesh, India \\ priyanshkhare0224@gmail.com}

\newtheorem{defi}{Definition}
\newtheorem{lemma}{Lemma}

\newtheorem{rem}{Remark}

\def\pf{\noindent{\it Proof.}}

\begin{document}

\maketitle

\fontsize{12}{15}\selectfont

\begin{abstract}
The article uses an Arithmetic-Geometric Fibonacci series to find the expected value, denoted \( E(\chi) \), of trials needed to observe \( k \) consecutive successes for the first time in a Bernoulli experiment using a recurrence relation. The article establishes that \( E(\chi) = 2(2^k - 1) \). It is important to note that this is not a new result, but to the best of my knowledge, this is a novel derivation of a well-established result. The other derivations of this result are cited in the references section.
\end{abstract}

\section{Introduction} Consider a coin-tossing experiment where the probability of success (i.e., getting heads) is 
$p$, and the probability of failure (i.e., getting tails) is $1-p$. For a fair coin, $p=\frac{1}{2}$.\\
This article explores the general case where $p$ is any real number within the interval $(0,1)$. We perform trials and record the outcome of each trial—heads or tails—sequentially. For instance, a possible outcome sequence in a $4$-trial experiment could be ``$hhth$,'' where $h$ represents heads and $t$ represents tails.
Our goal is to determine the average number of trials required to achieve $k$-consecutive successes for any positive integer $k$. For example, when $k=3$, a successful experiment may conclude in as few as $3$ trials with the outcome ``$hhh$.'' However, in other cases, it may take many more trials. For instance, the sequence ``$hhtththttth$" contains no occurrence of $3$-consecutive heads, so additional trials would be needed to achieve the first instance of $3$-consecutive heads.\\

\noindent A Bernoulli experiment consists of a series of independent Bernoulli trials, where each trial has a random outcome with two possible results: success or failure. The probability of success remains constant throughout the trials. Common examples of Bernoulli trials include tossing a coin, where the outcome is either heads or tails, or rolling a die, where one outcome is designated as success, while the others represent failure.\\

\noindent The aim of this article is to determine the expected number of trials (i.e., the average) needed to observe $k$-consecutive successes for the first time in a Bernoulli experiment. For the purposes of this paper, we refer to successes as ``heads" and failures as ``tails," analogous to a coin toss experiment. Using a generalized form of arithmetic, geometric, and Fibonacci series, we calculate the expected value of trials needed to observe $k$-consecutive successes for the first time.\\

\noindent While the results discussed here are not new and can be found in sources such as \cite{Kim}, \cite{Makri}, and \cite{Riggle}, this article is written by a high school student with the goal of making these concepts accessible to other high school students. It aims to demonstrate how basic ideas from arithmetic, geometric, and Fibonacci series, along with recurrence relations, can be applied to solve an intriguing problem: determining the average number of trials required to observe $k$-consecutive heads for the first time in a coin-tossing experiment.\\

\noindent In the next section, we discuss some preliminary concepts necessary to achieve our goal.
\section{Preliminaries}

\begin{defi}
For any positive integers $k$ and $n$, let $ \mathcal{F}_{n,k} $ be the set of sequences of heads and tails over $n$ trials, such that the first occurrence of $k$ consecutive heads appears in the $n$-th trial (i.e., the last $k$ symbols in the sequence are all heads). Denote the cardinality of $ \mathcal{F}_{n,k} $ by $|\mathcal{F}_{n,k} |=g_{n,k}$. 
\end{defi}
\begin{rem}\label{rem-initial-cond-fnk}
As implied by the above definition, we require $g_{n,k}=0$ whenever $n<k$ and $g_{n,k}=1$ whenever $n=k$. Furthermore, define $g_{0,k}=0$ for all positive integers $k$.
\end{rem}
\begin{lemma}
For any positive integers  $k$ and $n$ with $n\geq k$, the following holds:
\[
g_{n,k} = \sum_{m=1}^{k} g_{n-m,k}
\]
\end{lemma}

\pf 

To demonstrate that \( g_{n,k} = \sum_{m=1}^{k} g_{n-m,k} \), we will establish a bijection between the set \( \mathcal{F}_{n,k} \) and the union of the sets \( \{ \mathcal{F}_{n-i,k} \mid 1 \leq i \leq k \} \). Note that  by definition of $\mathcal{F}_{n,k}$,  $\mathcal{F}_{n-i,k}$ and $\mathcal{F}_{n-j,k}$ are disjoint for positive integers $i\leq k$ and $j\leq k$ whenever $i\neq j$.

Consider any valid sequence (i.e., an element) from \( \mathcal{F}_{n,k} \). In this sequence, let \( i \) be the position of the first appearance of tails. If the first tails appears at position \( i \), then the segment of the sequence followed by this tails has the length \( n-i \). This segment must be a valid sequence in \( \mathcal{F}_{n-i,k} \).

Now, consider a sequence in \( \mathcal{F}_{n-i,k} \). This sequence has \( k \) consecutive heads in the end by definition. These $k$ heads occupy positions beginning at the position \( n-i-k+1 \). By placing \( i \) tails in the beginning of this sequence, we obtain a sequence of length \( n \) where the first occurrence of \( k \) consecutive heads appears exactly at the $n$-th position. Thus, this extended sequence belongs to \( \mathcal{F}_{n,k} \). 

We have thus defined a bijection between \( \mathcal{F}_{n,k} \) and the union of the sets \( \{ \mathcal{F}_{n-i,k} \mid 1 \leq i \leq k \} \). Each sequence in \( \mathcal{F}_{n,k} \) uniquely corresponds to one in \( \mathcal{F}_{n-i,k} \) for some \( i \), and each sequence in \( \mathcal{F}_{n-i,k} \) can be uniquely extended to one in \( \mathcal{F}_{n,k} \) by concatenating a sequence of $i$ tails in the beginning.

As $\mathcal{F}_{n-i,k}$ and $\mathcal{F}_{n-j,k}$ are disjoint whenever $i\neq j$ and there is bijection between \( \mathcal{F}_{n,k} \) and  the set \( \{ \mathcal{F}_{n-i,k} \mid 1 \leq i \leq k \} \), the number of elements in \( \mathcal{F}_{n,k} \) equals the sum of the number of elements in the sets \( \mathcal{F}_{n-i,k} \) for \( i \) from 1 to \( k \):
\[
g_{n,k} = \left| \bigcup_{i=1}^{k} \mathcal{F}_{n-i,k} \right| = \sum_{i=1}^{k} |\mathcal{F}_{n-i,k}|
\]

Since \( |\mathcal{F}_{n-i,k}| = g_{n-i,k} \), it follows that:
\[
g_{n,k} = \sum_{m=1}^{k} g_{n-m,k}
\]
\qed

\section{Deriving an Expression for the Expectation Value}

For any positive integer $k$, consider performing a Bernoulli experiment, such as repeatedly tossing a coin, until we achieve $k$ consecutive successes for the first time. Let $S_k$ represent the set of all outcome sequences from this experiment. Specifically, these sequences end with exactly $k$ consecutive successes and contain no occurrence of $k$ consecutive successes prior to the final trial. Let \( \chi \) be the random variable that maps a sequence in $S_k$ to its length (i.e., to the number of trials in the sequence). Let $ p(\chi=i) $ denote the probability that a sequence in $S_k$ is of length $i$. For $i=n$ by definition of $F_{n,k}$,  for a fair coin with success probability in each trial $p=\frac{1}{2}$, each sequence in $F_{n,k}$ is equally likely. That is each sequence in $F_{n,k}$ is equally likely independent of number of heads or tails in the sequence. Hence, 

\begin{equation}\label{eq-1}
p(\chi=n) = \frac{|\mathcal{F}_{n,k}|}{2^n},
\end{equation}

The expectation value $ E(\chi)$ is expressed as:
\begin{eqnarray*}
E(\chi) & =& \sum_{i=1}^{\infty} i \cdot p(\chi = i)\\
\Rightarrow E(\chi) & = & \sum_{i=1}^{\infty} i \cdot \frac{g_{i,k}}{2^i}  \quad \mbox{[As $p (\chi=i) = \frac{|\mathcal{F}_{n,k}|}{2^n}$. ]}
\end{eqnarray*}
We now define the series \( y = \sum_{i=1}^{\infty} g_{i,k} \cdot r^i \). Note that the term-by-term differentiation of the series \( y \) with respect to \( r \) results in:
\[
\sum_{i=1}^{\infty} i \cdot g_{i,k} \cdot r^{i-1}
\]
\[
\Rightarrow E(\chi) = \frac{1}{2} \cdot y'
\]
To prove that the series \( y \) can be term-by-term differentiated, we must establish that it converges absolutely.

\subsection*{Proof of Convergence of the Series \( y \) Using the Ratio Test}

Consider the series \( y(r) = \sum_{i=k}^{\infty} g_{i,k} \cdot r^i \), with \( g_{i,k} \) defined recursively. Let \( a_i = g_{i,k} \cdot r^i \). Applying the ratio test, we compute:
\[
L = \lim_{i \to \infty} \left|\frac{a_{i+1}}{a_i}\right| = \lim_{i \to \infty} \left|\frac{g_{i+1,k} \cdot r^{i+1}}{g_{i,k} \cdot r^i}\right| = \lim_{i \to \infty} \left|\frac{g_{i+1,k}}{g_{i,k}} \cdot r\right|
\]
Substituting \( r = \frac{1}{2} \), we obtain:
\begin{eqnarray*}
L & = &\lim_{i \to \infty} \left|\frac{g_{i+1,k}}{2g_{i,k}} \right|\\
& = & \lim_{i \to \infty} \left|\frac{\sum_{j=1}^k g_{i+1-j,k} }{2g_{i,k}} \right|\\
& = & \lim_{i \to \infty} \left|\frac{\sum_{j=1}^k g_{i+1-j,k} }{g_{i,k}+g_{i,k}} \right|\\
& = & \lim_{i \to \infty} \left|\frac{\sum_{j=i-k+1}^i g_{j,k} }{g_{i,k}+g_{i,k}} \right|\\
&<& 1 \mbox{\quad [As $g_{i,k}=\sum_{j=i-k}^{i-1} g_{j,k}$ and $g_{t,k}>0$ for all $t\geq k$]}. 
\end{eqnarray*}
Thus, we conclusively establish \( L < 1 \) and prove that \( y \) converges absolutely via the ratio test.
\qed

\section{Deriving an Expression for \( y' \)}

Let \[
y = \sum_{i=1}^{\infty} g_{i,k} r^i
\]

where \( 0 < r < 1 \).  As $g_{i,k}=0$ for integers $i<k$, we get 
 \[
\Rightarrow y = \sum_{i=k}^{\infty} g_{i,k} r^i
\]

To find a solution, we multiply the equation by successive powers of \( r \) up to \( r^k \):

\[
yr = \sum_{i=k}^{\infty} g_{i,k} r^{i+1}
\]
\[
yr^2 = \sum_{i=k}^{\infty} g_{i,k} r^{i+2}
\]
\[
\vdots
\]
\[
yr^k = \sum_{i=k}^{\infty} g_{i,k} r^{i+k}
\]

Next, we subtract these equations from the original expression of \( y \). As $g_{i,k}=\sum_{j=1}^{k}g_{i-j, k}$ and $g_{i,k}=0$ for integers $i<k$, we get

\[
y \left( 1 - r - r^2 - \cdots - r^k \right) = r^k g_{k,k} 
\]
Note that all $g(i,k)-\sum_{j=1}^{k}g_{i-j, k}=0$ for all $i>k$ in the right side of the above equation. Evaluating the geometric series on the left side of the equation:
\begin{eqnarray*}
y \left( 1-\left( r+r^2+\cdots+r^k\right) \right) &=& r^k \mbox{\quad [As $g{k,k}=1$ for all positive integers $k$]}\\
\Rightarrow y \left(1-r\left(\frac{1-r^k}{1-r}\right)\right) & = & r^k\\
\Rightarrow y\left(\frac{1-r-r\left(1-r^k\right)}{1-r}\right) & = & r^k\\
\Rightarrow y\left(\frac{1-2r+r^{k+1}}{1-r}\right) & = & r^k\\
\Rightarrow y & = & \frac{r^k (1-r)}{1-2r+r^{k+1}}\\
\end{eqnarray*}
Therefore, 
\begin{eqnarray*}
\Rightarrow y' &= &\frac{\left(k r^{k-1}(1-r)-r^k\right)\left(1-2r+r^{k+1}\right)-r^k (1-r) \left(-2+(k+1)r^k\right)}{{\left(1-2r+r^{k+1}\right)}^2}\\
 &= &\frac{r^{k-1}\left(-r^{k+1}+2k r^2-3k r+k+r\right)}{{\left(1-2r+r^{k+1}\right)}^2}\\
\end{eqnarray*}
Finally, substituting the value of  \(E(\chi)\), we obtain:

\begin{eqnarray*}
 E(\chi) &=&  \frac{1}{2}\left( y'\right)\\
 & = &  \frac{1}{2}\left( \frac{r^{k-1}\left(-r^{k+1}+2k r^2-3k r+k+r\right)}{{\left(1-2r+r^{k+1}\right)}^2}\right)\\
\end{eqnarray*}
Below are some evaluations of expected values for $k$-run for various values of $k$ for a fair coin (i.e., $r=\frac{1}{2}$) based on the above obtained formula.
\begin{center}
     \begin{tabular}{| l | l | }
     \hline
     $k$ & $E(\chi)$  \\ \hline
     1 & 2  \\ \hline
     2 & 2  \\ \hline
     3 & 14  \\ \hline
     4 & 30 \\ \hline
      5 & 62 \\ \hline
     \end{tabular}
\end{center}
It can be easily shown using mathematical induction that for $r=\frac{1}{2}$, for any positive integer $k$ the value of $E(\chi)$ simplifies to:
\begin{eqnarray*}
 E(\chi) &=&  2(2^k-1).
\end{eqnarray*}.\qed
\section{Conclusion}
This paper derived the expectation value for the length of the first \( k \)-consecutive successes in a Bernoulli trial using arithmetic-geometric generalized Fibonacci series. By applying recurrence relations and some basic calculus, we established a formula for the expected number of trials required. Our final expression is:
\[
E(\chi) =  2(2^k-1).
\]
Deriving an analogous formula becomes significantly more challenging when the probability of success in a Bernoulli trial deviates from $\frac{1}{2}$ (i.e., when the coin is not fair), and no simplified closed-form solution exists. Note that for $r\neq \frac{1}{2}$, all sequences in the set $F_{i,k}$ are equally likely, which implies that $p(\chi=i)\neq \frac{g_{i,k}}{2^i}$ as given in equation (\ref{eq-1}). Besides $r\neq \frac{1}{2}$, another promising direction for further research could involve extending these findings to other stochastic processes.

\section{Acknowledgements}
I would like to extend my gratitude to Dr. Stoyan Dimitrov for introducing me to the problem, his advice, steering the direction and deciding the scope of this exploration. This paper would not have been possible without his support. I would also like to express my special thanks to Dr. Niraj Khare for putting me in touch with Dr. Stoyan Dimitrov and nurturing my love for mathematics. 

\end{document}